\documentclass[12pt]{amsart}
\usepackage{amsthm,amstext,amsmath,amscd,amssymb,latexsym,mathrsfs}

\usepackage{color}
\usepackage{amsfonts,array}
\usepackage{todonotes}
\usepackage{verbatim}
\usepackage{float}	
\usepackage{pgf,tikz}
\usepackage[toc,page]{appendix}
\usepackage[all,cmtip]{xy}
\usepackage{stmaryrd}
\usepackage{xparse}
\usepackage{listings}
\usepackage{enumitem}
\setlist{nosep, left=0pt}

\usepackage{tabularx,booktabs}

\usepackage{ragged2e} 
\usepackage{needspace}   
\usepackage{ltablex}   
\keepXColumns         
\usepackage{placeins} 

\newcolumntype{L}{>{\RaggedRight\arraybackslash}p{3.5cm}}
\newcolumntype{Y}{>{\RaggedRight\arraybackslash}X}

\usepackage{caption}
\usepackage{subcaption}
\usepackage[T1]{fontenc}
\usepackage{lmodern}
\usepackage{etoolbox}
\usepackage{xcolor}

\usepackage[shortalphabetic]{amsrefs}
\usepackage[all]{xy}
\usepackage{newtxtext}
\usepackage{newtxmath}
\usepackage{framed}
\usepackage{graphicx}
\usepackage{wasysym}
\usepackage{pdfpages}

\usepackage[thinlines,thiklines]{easybmat}

\usepackage{microtype} 

\usepackage{mathtools}  

\usepackage{bbm}
\usepackage[bookmarks=true, colorlinks=true, linkcolor=blue!50!black,
citecolor=GoetheBlue, urlcolor=orange!50!black, pdfencoding=unicode]{hyperref}
\usepackage[left=3.5cm,right=3.5cm,top=3cm,bottom=3cm,footskip=1cm
]{geometry}
\makeatletter
\def\thm@space@setup{%
  \thm@preskip=5pt \thm@postskip=5pt
}
\makeatother

\makeatletter
\renewcommand\@biblabel[1]{[#1]\quad}%

\@ifpackageloaded{amsrefs}{%
  \ifcsname BibLabel\endcsname
    \usepackage{etoolbox}
    \patchcmd{\BibLabel}{\hfill}{}{}{}%
    \renewcommand{\BibLabel}{%
      \ifcsname Hy@raisedlink\endcsname
        \Hy@raisedlink{\hyper@anchorstart{cite.\CurrentBib}\hyper@anchorend}%
        [\thebib]\quad
      \else
        [\thebib]\quad
      \fi
    }%
  \fi
}{}%
\makeatother

\newcommand{\restrict}{\,\rule[-5pt]{0.4pt}{12pt}\,{}}

\numberwithin{equation}{section}

\newcommand{\equ}{\ensuremath{\,=\,}}

\newcommand{\deq}{\ensuremath{\stackrel{\textrm{def}}{=}}}

\DeclareMathOperator{\pr}{pr}

\DeclareMathOperator{\lint}{int}

\newcommand{\Ra}{`\ensuremath{\Rightarrow}'  }
\newcommand{\La}{`\ensuremath{\Leftarrow}' }

\newcommand{\BA}{{\mathbb{A}}}

\newcommand{\BC}{{\mathbb{C}}}

\newcommand{\BF}{{\mathbb{F}}\,\!{}}
\newcommand{\BG}{{\mathbb{G}}}

\newcommand{\BN}{{\mathbb{N}}}

\newcommand{\BP}{{\mathbb{P}}}

\newcommand{\BR}{{\mathbb{R}}}

\newcommand{\BZ}{{\mathbb{Z}}}

\newcommand{\Fa}{{\mathfrak{a}}}
\newcommand{\Fb}{{\mathfrak{b}}}

\newcommand{\Fp}{{\mathfrak{p}}}
\newcommand{\Fq}{{\mathfrak{q}}}

\newcommand{\CC}{{\mathcal C}}

\newcommand{\CO}{{\mathcal O}}

\DeclareMathOperator{\Spec}{Spec}
\DeclareMathOperator{\Proj}{Proj}

\DeclareMathOperator{\supp}{supp}

\DeclareMathOperator{\im}{im}

\DeclareMathOperator{\Ann}{Ann}

\DeclareMathOperator{\rank}{rank}

\DeclareMathOperator{\Cone}{Cone}

\DeclareMathOperator{\Cox}{Cox}

\DeclareMathOperator{\Rel}{Rel}
\DeclareMathOperator{\Hom}{Hom}

\newcommand{\Quot}{\mathop{\rm Quot}\nolimits}

\newcommand{\Gen}{{\rm Gen}}

\definecolor{GoetheBlue}{RGB}{0,97,143}
\usepackage[bookmarks=true, colorlinks=true, linkcolor=blue!50!black,
citecolor=GoetheBlue, urlcolor=orange!50!black, pdfencoding=unicode]{hyperref}
\usepackage[left=3.5cm,right=3.5cm,top=3cm,bottom=3cm,footskip=1cm
]{geometry}
\makeatletter
\def\thm@space@setup{%
  \thm@preskip=5pt \thm@postskip=5pt
}
\makeatother

\newtheorem*{theorem*}{Theorem}

\newtheorem{theorem}{Theorem}[section]
\newtheorem{lemma}[theorem]{Lemma}

\newtheorem{proposition}[theorem]{Proposition}
\newtheorem{corollary}[theorem]{Corollary} 

\theoremstyle{definition}
\newtheorem{definition}[theorem]{Definition}
\newtheorem{example}[theorem]{Example}

\newtheorem{remark}[theorem]{Remark}

\theoremstyle{plain}
\newtheorem{introthm}{Theorem}

\def\phi{\varphi}
\def\epsilon{\varepsilon}
\def\setminus{\smallsetminus}
\let\oldbullet\bullet
\def\bullet{{\mathchoice{\oldbullet}%
                        {\oldbullet}%
                        {\scriptscriptstyle\oldbullet}%
                        {\oldbullet}}}
\let\oldemptyset\emptyset
\let\emptyset\varnothing

\title{Torsor and Quotient Presentations for $D$-homogeneous Spectra} 

\author{Felix G\"obler}
\address{Institut f\"ur Mathematik, Goethe-Universit\"at Frankfurt, Robert-Mayer-Str. 6-10, D-60325 Frankfurt am Main, Germany }
\email{goebler@math.uni-frankfurt.de}

\begin{document}

\pagenumbering{arabic} 
\setcounter{page}{1}

\maketitle





\begin{abstract}
The $D$-graded Proj construction provides a general framework for constructing schemes from rings graded by finitely generated abelian groups $D$, yet its properties and applications remain underdeveloped compared to the classical $\BN$-graded case. This paper provides the essential characteristics of $D$-graded rings $S$, like the distinction between $D$-homogeneous prime ideals and $D$-prime ideals if $D$ has torsion. 

We particularly focus on describing the quotient by the associated group scheme, generalizing the construction of a toric variety from its Cox ring.
As in the $\BN$-graded construction, the basic affine opens of the Proj construction are given in terms of degree-zero localizations $S_{(f)}$, where $f$ in $S$ homogeneous is \emph{relevant}.
We prove that  $\pi_f: \Spec(S_f) \to \Spec(S_{(f)})$ is a geometric quotient if $f$ is relevant, and give necessary and sufficient conditions for this map to be a pseudo $\Spec(S_0[D])$-torsor. 
\end{abstract}

\tableofcontents


\section*{Introduction}
Projective space is the basic example of a variety obtained from a graded ring, and the usual \(\Proj\) construction for \(\BN\)-graded rings is the standard algebraic tool that realizes this correspondence. 
However, many natural geometric objects carry rings graded by finitely generated abelian groups (called \emph{multigraded rings}) rather than by \(\BN\).
Concretely, this is the case for Cox rings of toric (pre)varieties, and more generally, coordinate rings of divisorial varieties or Mori dream spaces.
These multigradings encode divisor class information and symmetries that cannot be captured by the classical \(\BN\)-graded Proj.

Despite the ubiquity of multigraded rings in modern algebraic geometry, their projective geometry remains underdeveloped compared to the classical theory. The $D$-graded $\Proj$ construction, introduced by \cite{BS}, provides the natural framework for studying these objects, yet its foundational aspects and geometric interpretation have received limited attention.

This paper aims to fill this gap by providing a comprehensive treatment of the properties of $D$-graded rings, with particular focus on the quotient description of the corresponding $D$-graded Proj construction.
Specifically, that quotient description recovers familiar quotient and coordinate descriptions (in particular, Cox-type quotients).
The geometric and topological interpretation of the $D$-graded Proj construction, as well as the strong connection to toric prevarieties, are discussed in \cite{paper2}, \cite{paper4}, and \cite{paper3} (which also contains an overview of the research associated with the Brenner--Schröer Proj construction).

Let $D$ be a finitely generated abelian group and
\[
S \equ \bigoplus_{d\in D} S_d
\]
a commutative unital $D$-graded ring. Following \cite{BS}, the homogeneous elements $f \in S$, whose group of units $D^f$ in $S_f$ has finite index in $D$  are called \emph{relevant}. Then the $D$-homogeneous spectrum is defined by
\[
\Proj^D(S) \deq \bigcup_{f \in S \text{ relevant}}\Spec\bigl(S_{(f)}\bigr),
\]
and the \emph{irrelevant ideal} $S_+$ is defined to be the ideal generated by all relevant elements in $S$.
The relevance of $f$ ensures that the basic affine opens $\Spec(S_{(f)})$ are of equal dimension $\dim(S) - \rank(D)$ (see Lemma~\ref{lem:dim_proj}).


$\Proj^D(S)$ is not just a formal generalization of ordinary Proj, it is the correct coordinate ring picture for many important non-affine varieties. In the toric case, for example, every simplicial toric variety $X$ is naturally isomorphic to the geometric quotient $(\Spec(\Cox(X)) \setminus Z) \sslash G$ (see \cite{Cox}, Theorem 2.1). The corresponding generalization for $D$-graded rings is the following.

\begin{introthm}[{\autoref{lem:periodic_quotient}}]
Let $S$ be an effectively $D$-graded ring.
Then for all relevant $f \in S$, the map $\pi_f \colon \Spec(S_f) \to \Spec(S_{(f)})$ induced by $S_{(f)} \to S_f$ is a geometric quotient in the category of schemes (in the sense of \cite[\href{https://stacks.math.columbia.edu/tag/04AE}{Definition 04AE}]{stacks-project}). 
In particular, these patches glue to the geometric quotient
\begin{align*}
    \Proj^D(S) \stackrel{\sim}{\equ} D(S_+) / \Spec(S_0[D])
\end{align*}
\end{introthm}

Hence, the GIT presentation of simplicial toric varieties is indeed a special case of the $D$-graded Proj construction.

In addition, we examine the related concept of a \emph{pseudo $G$-torsor} (see  \cite[\href{https://stacks.math.columbia.edu/tag/0498}{Definition 0498}]{stacks-project}),  also studied in \cite{May}, Proposition 5.13. We prove a stronger result.

\begin{introthm}[\autoref{thm:quotient_pi_f_general}]
    Let $S$ be effectively $D$-graded, $f \in S$ relevant and $\pi_f\colon \Spec(S_f) \to \Spec(S_{(f)})$ be the map induced by $S_{(f)} \to S_f$. Then $\Spec(S_f)$ is a pseudo $G$-torsor if and only if $\deg( (S_f)^\times) = D$ (or equivalently $D^f = D$) and $S$ is integral.
    In particular, if $\deg( (S_f)^\times) = D$ holds and $S$ is integral, then $\pi_f$ is Zariski locally trivial for the base change $S_0 \to S_{(f)}$.
\end{introthm}

{\normalfont\bfseries Acknowledgements.}
First, I would like to thank my supervisor, Alex Küronya, for his tremendous support and insightful feedback on this version.
I am also very grateful to Johannes Horn, Andrés Jaramillo Puentes, Kevin Kühn, Jakob Stix, Martin Ulirsch, and Stefano Urbinati for many helpful discussions and valuable comments.
In particular, I want to thank an anonymous referee for his suggestions.

Partially funded by the Deutsche Forschungsgemeinschaft (DFG, German Research Foundation) TRR 326 \textit{Geometry and Arithmetic of Uniformized Structures}, project number 444845124, and by the LOEWE grant \emph{Uniformized Structures in Algebra and Geometry}.

\section{$D$-graded Proj}\label{sec:mult_ring}
This paper aims to give an overview of the foundations and abstract properties of the $D$-graded Proj construction. The starting point is the following datum.

\begin{definition}\label{def:graded_ring}       
Let $D$ be a finitely generated abelian group. 
\begin{enumerate}[label=(\arabic*)]
    \item A \emph{$D$-graded} ring is a commutative ring $S$ with a decomposition $S = \oplus_{d \in D} S_d$ into a direct sum of subgroups $S_d \le S$ such that $S_d \cdot S_e \subseteq S_{d+e}$. 
The elements of $S_d$ are called \emph{$D$-homogeneous of degree $d$}. In particular, $0 \in S$ is homogeneous of every degree.
We will call elements \emph{homogeneous} instead of $D$-homogeneous when there is no ambiguity and denote the set of $D$-homogeneous elements by $h_D(S)$, respectively $h(S)$.

    \item A commutative ring $S$ with a grading by a finitely generated abelian group is called a \emph{multigraded ring}\footnote{If we fix a finitely generated abelian group $D$, we speak of a \emph{$D$-graded ring} to emphasize that the grading with respect to $D$ is part of the data. In other words, every $D$-graded ring is multigraded, but when we specify $D$, we regard the grading itself as explicit and fixed.}.
    
    \item Let $R$ be $D_R$-graded and $S$ be $D_S$-graded rings, where $D_R$ and $D_S$ are finitely generated abelian groups. A \emph{morphism of multigraded rings} is a pair $(\varphi, \alpha_\varphi)$, where $\varphi$ is a ring homomorphism $\varphi\colon R\to S$, and $\alpha_\varphi$ is a group homomorphism $\alpha_\varphi\colon D_R \to D_S$ satisfying $\varphi(R_d) \ \subseteq \ S_{\alpha_\varphi(d)}$ for all $d\in D_R$. 
\end{enumerate}
Formally, a multigraded ring is a pair $(D, S)$, where $D$ is a finitely generated abelian group and $S$ is a $D$-graded commutative unital ring.
For readability, we will usually suppress the grading group and denote a multigraded ring by its underlying ring $S$.
\end{definition}

Since $S = \oplus_{d\in D}S_d$, every element $f \in S$ has a unique decomposition $f = \sum_{d\in D} f_d$ into the \emph{$D$-homogeneous components} $f_d$ of $f$ of degree $d$, where $f_d = 0$ for almost all $d\in D$.
We say that an ideal $\Fa \unlhd S$ is \emph{$D$-graded}, if $f = \sum_{d \in D} f_d \in \Fa$ implies that $f_d \in \Fa$ for all $d\in D$. If the group $D$ is clear from the context, we might just say that $\Fa$ is graded.

\subsection{Multigraded Rings}\label{subsec:1.1}

In this section, we analyze some of the basic properties of multigraded rings. Note that those details have not been addressed in the literature yet. We closely follow \cite{Bosch}, 9.1/2 up to 9.1/6, which describes the $\BN$-graded case that we aim to generalize.

\begin{lemma}\label{lem:bosch_9.1/2}
An ideal $\Fa \unlhd S$ is $D$-graded if and only if it is generated by $D$-homogeneous elements.
As a direct consequence, sums, products, and intersections of $D$-graded ideals remain $D$-graded ideals. 
\end{lemma}

\begin{proof}
The first statement holds by  the same arguments as found in \cite{Bosch}, Remark 9.1/2.
Hence, also sums and products remain $D$-graded.
As the intersection of two $D$-graded ideals $\Fa$ and $\Fb$ also contains all homogeneous components of elements of $\Fa$ respectively $\Fb$, the claim follows.
\end{proof}

However, the radical of a $D$-graded ideal is not $D$-graded in general:

\begin{example}\label{ex:radical_not_D_hom}
Let $p$ be a prime number and let $S := \BF_p[x]/(x^p-1)$ be graded by $D := \BZ/ p\BZ$ such that $\deg(x) = \overline{1}$.
Then $(0)$ is a homogeneous ideal and
\begin{align*}
    (x-1)^p \equ x^p -1 \equ 0,
\end{align*}
that is, $x-1 \in \sqrt{(0)}$. But neither $x$ nor $1$ are nilpotent, i.e. $1, x \not\in\sqrt{(0)}$, therefore $\sqrt{(0)}$ cannot be $D$-graded.
\end{example}

\begin{remark}
 Let $\Fa \subset S$ be a $D$-homogeneous ideal. 
 If $D$ admits a total order compatible with the group structure, then $\sqrt{\Fa}$ remains $D$-homogeneous. In particular, this rule applies for $D=\BZ^r$ together with the lexicographic order.
\end{remark}

We want to consider ideals of the following type (cf.\ \cite{BS}, Remark 2.3, and \cite{CRB}, Definition 1.5.3.1 (iv)).

\begin{definition}\label{def:D_prime}
Let $\Fp \unlhd S$ be a proper graded ideal and $h_D(S)$ denote the set of $D$-homogeneous elements of $S$. Then $\Fp$ is called a \emph{$D$-prime ideal} if, for $D$-homogeneous elements $f,g\in h_D(S)$ with $fg\in\Fp$, either $f\in\Fp$ or $g\in\Fp$.
\end{definition}

In the $\BN$-graded case we know that, if $\deg_D(S) \subseteq \BZ$ and $\Fp \unlhd S$ is a proper graded ideal, then $\Fp$ is $D$-prime if and only if $\Fp$ is prime (this is \cite{Bosch}, 9.1/4). 
In the $\BZ^r$-graded case, we can maintain this characterization of graded prime ideals, as shown in \cite{Rob}, Lemma 2.1.1.

\begin{proposition}\label{prop:D_prime_prime_Z^r}
Let $D = \BZ^r$, let $S$ be a factorially $\BZ^r$-graded\footnote{Factorially graded means that homogeneous elements have unique homogeneous factorization.} integral noetherian domain, where $r \le \dim(S)$\footnote{If $r > \dim(S)$, $S$ has no relevant elements and $\Proj^D(S)$ is empty.}, and let $\Fp \unlhd S$ be a proper graded ideal. Then $\Fp$ is a prime ideal if and only if $\Fp$ is a $D$-prime ideal.
\end{proposition}

In general, only $\Ra$ holds, as the following example for finite $D$ shows.

\begin{example}\label{ex:D_prime_not_prime}
Let $S = \BR[x]$, let $D = \BZ/2\BZ$ and $\Fp = (x^2-1)$. Then $S = S_{\overline{0}} \oplus S_{\overline{1}}$, where $S_{\overline{i}} = \sum_{n \ge 0} \BR x^{2n+i}$. The ideal $\Fp$ is not prime, since $(x+1)(x-1) = 0$ in $S/\Fp$ but $(x\pm 1) \neq 0$ in $S/\Fp$. On the other hand, $x \pm 1$ is not $D$-homogeneous and $\Fp$ is graded since $x^2-1 \in S_{\overline{0}}$. 
In particular, $\Fp$ is $D$-prime: If $fg \in\Fp$ for homogeneous $f, g \in S$ then $(x^2-1)\mid fg$ and it follows $f \in \Fp$ or $g \in \Fp$ by degree arguments. On the other hand, $\Fq = (x^2+1)$ is a homogeneous prime ideal of $S$. 
\end{example}

As in the $\BN$-graded case, the $D$-graded Proj is defined in terms of degree-zero localizations. 

\begin{definition}
    For $f \in S_d$, the ring $S_f$ is defined to be the graded ring generated by the subgroups
\begin{align*}
    (S_f)_e \deq \left\{\frac{s}{f^k} \mid k \in \BN, s \in S_{e+kd} \right\}, \quad e \in D,
\end{align*}
in analogy with \cite{Bosch}, Proposition 9.1/5.
The subring $S_{(f)} := (S_{f})_0 \subset S_f$ is called the \emph{$D$-homogeneous localization} of $S$ by $f$. Likewise, we define the $D$-homogeneous localization of a $D$-prime ideal $\Fp \unlhd S$ to be
\begin{align*}
    S_{(\Fp)} \deq \left\{\frac{g}{f^k} \mid k \in \BN, f \not\in \Fp \text{ relevant}, g \in S_{k \cdot \deg(f)}\right\}.
\end{align*}
\end{definition}

In the classical setting, two homogeneous elements $f, g \in S$ are always `related', in the sense that there are integers $k, l$ such that $\deg(f^k) = \deg(g^l)$. It is immediate that the degree-zero localizations highly depend on the existence of such equations. In a multigraded ring, there can exist homogeneous elements that are not (integral) linearly dependent.

\begin{example}\label{ex:double_origin_P^1_first}
    Let $S = \BC[x, y, z]$ graded by $\BZ^2$, where $x \mapsto (1, 0)$, $y \mapsto (0, 1)$ and $z \mapsto (1, 1)$. Then the degrees of $x$ and $y$ are linearly independent in $D_\BR = \BR^2$. We claim that $S_+ = (xy, xz, yz)$ (we will prove this in general in Lemma~\ref{lem:mon_gen_S+}) and compute the corresponding degree-zero localizations.
    It holds that
    \begin{align*}
        S_{(xy)} \equ \BC[\frac{z}{xy}],\ S_{(xz)} \equ \BC[\frac{xy}{z}] \text{ and } S_{(yz)} \equ \BC[\frac{xy}{z}].
    \end{align*}
    In particular, we see that the generators of the degree-zero localizations correspond to the linear dependencies between the degrees of the variables in $D_\BR$, i.e.\ $\deg(xy) = \deg(z)$. It also shows why we do not want to take arbitrary homogeneous elements for degree-zero localizations, as $S_{(f)} = \BC$ for $f= x, y$.
\end{example}


\subsection{Relevant Elements}\label{sec:relevant_elements}
In the classical $\BN$-graded Proj construction, the irrelevant ideal is defined as the positively graded part. As in general, $D$ does not need to be totally ordered, we may not have a notion of positivity. Thus, Brenner--Schröer defined the irrelevant ideal $S_+$ to be generated by all relevant elements, i.e.\ homogeneous elements $f$ such that the units in $S_f$ form a finite index subgroup in $D$. 

One may think of relevant elements as elements having sufficiently many homogeneous divisors, as those divisors exactly correspond to units in $S_f$ (cf.\  Lemma~\ref{lem:relevant}). We also present geometric criteria in terms of the \emph{weight cones} $\CC_D(f)$ in $D_\BR$, that are generated by the degrees of homogeneous divisors of (powers of) $f$.

Suppose $S = R[T_1,\ldots, T_n]$, $\rank(D) = r < n$, is a polynomial ring, where $R$ is a commutative unital ring. In that case, the minimal relevant generators of the irrelevant ideal (i.e.\ the ideal generated by all relevant elements) are monomials $f$, given by $f = \prod_{i=1}^r T_{j_i}^{\epsilon_i}$, $j_i \in \{1, \ldots, n\}$, such that the $\deg(T_{j_i})$ are pairwise linearly independent (cf.\ Lemma~\ref{lem:mon_gen_S+}), where the $\epsilon_i$ are either all equal to 1 or all equal to zero (if $1$ is relevant).

Most importantly, relevant elements $f$ give affine toric varieties if $S$ is a noetherian polynomial ring, i.e.\ $\Spec(S_{(f)})$ is a simplicial affine toric variety (see \cite{BS}, Proposition 3.4). 
The following definition is due to \cite{BS}. If not explicitly stated otherwise, $D$ will always denote a finitely generated abelian group and $r = \rank(D)$.

\begin{definition} \label{def:relevant}
Let $S$ be a $D$-graded ring. 
\begin{enumerate}[label=(\arabic*)]
    \item The ring $S$ is called \emph{periodic} if the degrees of the homogeneous units $f \in S^\times \cap h(S)$ form a subgroup $D'\le D$ of finite index, that is
    \begin{align*}
        D' \deq \{\deg_D(f) \mid f \in S^\times \text{ homogeneous }\}\le D
    \end{align*}
     is of finite index.
    \item An element $f \in S$ is called \emph{relevant} if
    \begin{enumerate}[label=(\roman*)]
        \item $f$ is homogeneous and
        \item the localization $S_f$ is \emph{periodic}.
    \end{enumerate}
    We denote the set of relevant elements of $S$ by $\Rel^D(S)$.
    \item Given a homogeneous element $f \in S_d$, we define the \emph{support of $f$} to be the subset of $D$ given by
\begin{align*}
    C^f \deq \{\deg_D(g) \mid g \text{ is a homogeneous divisor of $f^k$ for some $k\ge 0$}\}.
\end{align*}
Furthermore, we define the \emph{weight cone of $f$} to be the closed convex cone
\begin{align*}
    {\CC_D(f)} &\deq \overline{\Cone}(C^f) 
\end{align*}
in $D_\BR := D \otimes_\BZ \BR$. If $f$ is relevant, we also call $\CC_D(f)$ a \emph{relevant weight cone}.

\item For homogeneous $f \in h(S)$, we define the \emph{support group} of homogeneous units in $S_f$ to be
\begin{align*}
    D^f \deq \{\deg_D(g) \mid g \in (S_f)^\times \cap h(S_f)\}.
\end{align*}
\item Finally, we define the \emph{irrelevant ideal $S_+$} of $S$ to be the ideal generated by all relevant elements\footnote{The chosen language seems counterintuitive, but it is nevertheless the generalization of the $\BN$-graded case. 
The designation results from the interpretation of $V(S_+)$ as `irrelevant locus' and of $D(S_+)$ as `relevant locus', and from the fact that $\Proj^D(S)$ coincides with the quotient of $D(S_+)$ by $\Spec(S_0[D])$.} of $S$, i.e.\ $S_+ = (\Rel^D(S))$. 
\end{enumerate} 
\end{definition}

\begin{remark}
    Let $f \in S$ be homogeneous. 
    \begin{enumerate}
        \item 
     By definition, $S_f$ has a unit in every degree $d \in D^f$. In particular, $f$ is relevant if and only if for all $d \in D$, $(D:D^f) d$ is the degree of a unit in $S_f$.
        \item The support $C^f$ of $f$ is a monoid, and the support group $D^f$ is indeed a group. In particular, $D^f$ is exactly the group generated by the monoid $C^f$. 
    \end{enumerate}   
\end{remark}

Relevance depends on the context in which the grading is understood.

\begin{example}\label{ex:effective_needed}
    Let $S = \BC[x, y]$ and $D = \BZ$ with $\deg(x) = 1$ and $\deg(y) = 1$. Then $x$ and $y$ are relevant. However, we might want to view $D$ embedded in $\BZ^2$, such that $\deg(x) = \deg(y) = (1, 0)$. In this case, $x$ and $y$ are no longer relevant. But the actual support of the two gradings does not differ.
\end{example}

Therefore, we restrict to \emph{effective} gradings (cf.\  \cite{CRB}, Definition 1.1.1.4), so that the grading group is forced to be generated by its support.

\begin{definition}\label{def:grading_effective}
    Let $S$ be a $D$-graded ring, where $D$ is a finitely generated abelian group. 
\begin{enumerate}[label=(\arabic*)]
    \item The \emph{weight monoid} of $S$ is defined to be
    \begin{align*}
        \omega(S) \deq \{ d \in D \mid S_d \neq 0\} \subseteq D .
    \end{align*}
    \item The $D$-grading of $S$ is called \emph{effective}, if the weight monoid $\omega(S)$ generates $D$ as a group.
    \item The \emph{weight cone} of $S$ is defined to be
    \begin{align*}
        \sigma(S) \deq \overline{\Cone}(\omega(S)) \subseteq D \otimes_\BZ \BR
    \end{align*}
\end{enumerate}
\end{definition}

From now on, we will always assume that the $D$-grading is effective. Otherwise, we may replace $D$ by the group generated by $\omega(S)$.

The following characterizations are due to \cite{BS}, \cite{KU}, and \cite{May} (mostly without proof). Therefore, we give a proof whenever there is none in the literature, and cite otherwise.

\begin{lemma}[Criteria for relevance]\label{lem:relevant} 
Let $f \in S$, $f \neq 1$, be homogeneous. Then the following are equivalent:
\begin{enumerate}[label=(\arabic*)]
    \item $f$ is relevant.
    \item  The degrees of all homogeneous divisors $g \mid f^n$, $n\ge 0$, \emph{generate} a subgroup $D'\le D$ of finite index (in fact, $D'= D^f)$.
    \item For all $d \in D$ there exists a positive integer $N > 0$ such that $Nd \in D^f$.
    \item $\lint(\CC_D(f)) \subseteq D_\BR$ is non-empty.
    \item $\deg(f) \in \lint(\CC_D(f))$. 
\end{enumerate}
Furthermore, relevance is closed under multiplication and also closed under multiplication with $S_0$.
\end{lemma}

\begin{proof}
The equivalence of $(1)$ and $(2)$ is straightforward, as every unit in $S_f$ corresponds to a divisor of (some power of) $f$ and vice versa.

$(3)$ is basically a reformulation of $(1)$, where $N$ is the index of $D^f$ in $D$.

Regarding $(4)$, one can see that the maximality of $\dim(\CC_D(f))$ is equivalent to $\rank(D) = \rank(D^f)$, which in turn is equivalent to $D^f$ having finite index in $D$.

For $(5)$, we use \cite{KU}, Lemma 2.10 for \Ra. 
Conversely, if $\deg(f) \in \lint(\CC(f))$, then $\lint(\CC(f)) \neq \oldemptyset$ and $f$ is relevant by (3).

Now, using $(4)$, one deduces that relevance is closed under multiplication, as for $f_1, f_2 \in S$ relevant, it holds $\CC_D(f_i) \subset \CC_D(f_1f_2)$. The same argument shows closedness with respect to multiplication by $S_0$.
\end{proof}

\begin{example}\label{ex:double_origin_P^1_2nd}
    Consider $S = \BC[x, y, z]$ graded by $\BZ^2$ with $x \mapsto (1, 0)$, $y \mapsto (0, 1)$ and $z \mapsto (1, 1)$ from Example~\ref{ex:double_origin_P^1_first}, where we claimed that $S_+ = (xy, xz, yz)$. Note that for $f = x, y, z$, $D^f$ is generated by the units $f$ and $f^{-1}$, so the units in $S_f$ do not give a finite index subgroup in $D$.
    It is also immediate that the cones $\CC_D(f)$ for those $f$ are one-dimensional.
    Therefore, those elements cannot be relevant. So the next step is to take all monomials that are given by the product of two variables, i.e.\ $f = xy, xz, yz$. It is not hard to see that every localization $S_f$ contains enough units. We will show this for $f = xy, xz$: \\
    For $f = xy$, $S_f$  has units $x, y, \frac{1}{x}, \frac{1}{y}$ in degree $\pm e_i$ for $i = 1,2$. Thus $S_f$ contains a unit in every degree $d \in \BZ^2$, i.e.\ $D^f = D = \BZ^2$. \\
    For $f = xz$, there are units $x, z, \frac{1}{x}, \frac{1}{z}$ in degree $\pm e_1$ and $\pm (e_1+e_2)$ in $S_f$. In particular, $\frac{z}{x}$ is a unit of degree $e_2$, so that again $D^f = D$ and $S_f$ contains a unit in every degree $d \in D$.
    
    Now, as $x, y, z$ are not relevant, but $xy, xz, yz$ are, we deduce that $S_+ = (xy, xz, yz)$.
\end{example}

Lemma~\ref{lem:relevant} (3) has an important consequence:

\begin{corollary}\label{cor:rel_many_el_deg_zero}
    Let $h \in h(S)$ be homogeneous and $f \in S$ relevant. Then there exists an element $g_h \in h(S)$, an element $k \in \BZ$ and an element $N > 0$ such that
    \begin{align*}
        \deg\left(\frac{h^N g_h}{f^k}\right) \equ 0 .
    \end{align*}
\end{corollary}

\begin{proof}
    By Lemma~\ref{lem:relevant} (3), there exists an $N > 0$ and a homogeneous $f' =\frac{a}{f^l} \in (S_f)^\times$, such that $N\deg(h) = \deg(f')$, where $a \in S$ is homogeneous by definition and $\deg(f') = \deg(a) - l\deg(f)$. 
    Since $f'$ is a unit in $S_f$, there exists a homogeneous $g_h\in S$ and an $m\in\mathbb N$ such that $ag_h=f^m$. Hence
\[
\deg(g_h)=m\deg(f)-\deg(a).
\]
Using $N\deg(h)=\deg(a)-l\deg(f)$, we get
\[
N\deg(h)+\deg(g_h)=(m-l)\deg(f).
\]
Thus
\[
\deg\left(\frac{h^N g_h}{f^{m-l}}\right)=0,
\]
and the claim follows.
\end{proof}

We want to choose generators of $S_+$ when $S$ is a noetherian polynomial ring or a noetherian factorially graded domain (see \cite{CRB}, Definition 1.5.3.1).

\begin{definition}\label{def:length_fact}
Let $0 \neq f \in S = R[T_1,\ldots, T_n]$ be homogeneous and $k \in \BN$ such that $k \le r$, where $R$ is a unital commutative ring and $r = \rank(D)$.
\begin{enumerate}
    \item We say that $f$ has a \emph{length $k$ factorization} if there is a subset $I \subseteq \{1, \ldots, n\}$ of size $k$, an element $u \in R$ and integers $\epsilon_i \in \BN_+$, $i \in I$, such that 
    \begin{align*}
        f \equ u \cdot \prod_{i\in I} T_i^{\epsilon_i}
    \end{align*}
    where the $d_i = \deg(T_i)$ generate a subgroup of $D$ of rank $k$.

        \item We say that $f$ has a \emph{linear length $k$ factorization} if there is a subset $I \subseteq \{1, \ldots, n\}$ of size $k$, and an element $u \in R$  such that 
    \begin{align*}
        f \equ u \cdot \prod_{i\in I} T_i,
    \end{align*}
    where the $d_i$, $i\in I$, generate a subgroup of rank $k$.

    \item Let $f$ be relevant. We say that $f$ is \emph{monomic} if $f$ has a linear length $r$ factorization and $u = 1_R$.
\end{enumerate}
\end{definition}

\begin{example}
    Let $S = \BC[x, y, z]$ be graded by $\BZ^2$, where $x \mapsto e_1$, $y \mapsto e_2$ and $z \mapsto e_1 + e_2$. Then the elements $xy, xz, yz \in S_+$ are monomic and minimal, i.e. there is no relevant element $g$ of smaller degree being distinct to $xy, xz, yz$ that divides $xy, xz$ or $yz$.
    The element $f = x(xy+z)$, for example, is relevant, but not monomic. 
Note that $f = x^2y + xz$ is a sum of homogeneous elements of the same degree, where $f \in (xy, xz, yz)$. 
\end{example}

We prove the following characterization of relevant elements for polynomial rings. It can be seen as the generalization of \cite{CRB}, Proposition 2.1.3.4.

\begin{lemma}\label{lem:mon_gen_S+}
Let $R$ be a UFD and $S =R[T_1, \ldots, T_n]$ be a $D$-graded polynomial ring, such that $1$ is not relevant and $\rank(D) = r < n$.
    Then the irrelevant ideal $S_+$ is generated by all monomic relevant elements, and there are at most $\binom{n}{r}$ monomic relevant elements.
\end{lemma}

\begin{remark}\label{rem:monomic_dep}
Fix a homogeneous presentation \(S=R[T_1,\dots,T_n]\) with \(\deg T_i=d_i\in D\). 
The description of \(S_+\) in Lemma \ref{lem:mon_gen_S+} is relative to this presentation:
The ideal \(S_+\) (and hence \(\operatorname{Proj}^D(S)\)) is intrinsic to the graded ring \((D,S)\).
But the specific monomial generating set need not be preserved by arbitrary (non-graded) $R$-algebra automorphisms. It is canonical only up to graded automorphisms of \((D, S)\).
But the same holds for the Proj construction itself, so we do not have to worry about non-graded automorphisms at all.
\end{remark}

\begin{proof}
    First, note that $S$ is factorially $D$-graded, since $R$ is a UFD.
    Let $f \in S$ be relevant. Without loss of generality, we may assume that $f = h_1 \ldots h_r$ for irreducible homogeneous elements $h_1, \ldots, h_r$ (if $f$ has a factorization with more than $r$ irreducible factors, the proof goes through verbatim). If all $h_i$ are monomials, there is nothing to show. Thus assume  that some $h_j$, $j \in \{1, \ldots, r\}$, is not a monomial. Then $h_j$ has to be a sum of homogeneous elements of equal degree, so we may assume there are indices $i, k, l \in \BN$ such that $h_j$ is of type $h_j = (T_k T_l - T_j)$, where $\deg(T_kT_l) = \deg(T_j) = \deg(h_j)$ (any more complex equation like $(T_k T_l - T_jT_i)$ will result in the same argument).
    Note that the precise form of the relation is determined by the grading of $S$, but it must be a sum of monomials of same degree.
    Hence, the argument does not really depend on the form of the equation - we might just get more than $r$ factors - which is fine.
    Since $f = h_1 \ldots h_r$ and the degrees of the $h_i$ are pairwise linearly independent, every homogeneous summand of $f$ has a factorization of length $\ge r$ (after expanding). 
    More concretely, each summand arises by choosing one homogeneous summand from each of the $r$ linearly independent factors. Hence, every summand is a product of at least $r$ homogeneous (indeed irreducible) pieces and therefore an element of the ideal generated by the monomic relevant elements.
    The claim follows.
    Furthermore, $S_+$ is a radical ideal, as monomic relevant elements are square-free by definition. 
\end{proof}

As every relevant element is given by a finite linear combination of monomic relevant elements by Lemma~\ref{lem:mon_gen_S+}, monomic relevant elements are minimal generators for $S_+$.
We will denote a choice of a minimal generating set by monomic relevant elements for $S_+$ by $\Gen^D(S)$.

\begin{remark}
    Lemma~\ref{lem:mon_gen_S+} has a direct generalization to factorially graded noetherian rings. 
\end{remark}

\begin{example}\label{ex:standard_ex}
We want to use the previous Lemma to compute monomic relevant generators of $S_+$. 
\begin{enumerate}
    \item Let $S = \BC[x, y, z, w]$ and $D = \BZ^2$, where $\deg(x) = \deg(y) = (1, 0)$, $\deg(z) = (1, 1)$ and $\deg(w) = (0, 1)$. It holds that $\Gen^D(S) = \{xw, yw, zw, xz, yz\}$. 

    \item Let $S = \BC[x, y, z]$ be graded by $\BZ\times \BZ/2\BZ$, such that $\deg(x) = (1, 0)$, $\deg(y) = (0, 1)$ and $\deg(z) = (1, 1)$. Then $y$ is no longer relevant (compared to Example~\ref{ex:double_origin_P^1_first}), as $D^y = \BZ/2\BZ$ cannot have finite index in $D$. In particular, monomic relevant elements have length one. 
    Thus, $\Gen^D(S) = \{x, z\}$. 
\end{enumerate}    
\end{example}


\subsection{$D$-graded Proj Construction}\label{sec:BS_Proj}

Having dealt with the notion of relevance, we can introduce the $D$-graded Proj construction.
A $D$-grading on $S$ gives rise to a coaction
\begin{align*}
    S \to S \otimes_{S_0} S_0[D], \ \sum_{d \in D} s_d \mapsto \sum_{d\in D} s_d \otimes d.
\end{align*}

If we view the group ring $S_0[D]$ as free $S_0$-module with basis $D$, and define the multiplication of $\alpha, \beta \in S_0[D]$ to be
\begin{align*}
    (\sum_i a_i d_i) \cdot (\sum_j b_j e_j) \deq \sum_{i,j} a_i b_j d_i e_j ,
\end{align*}
$S_0[D]$ becomes an $S_0$-algebra, called \emph{group algebra of $D$}.

By defining $\Delta(\alpha) = \alpha \otimes \alpha$ (\emph{comultiplication}), $\epsilon(\alpha) = 1$ (\emph{counit}) and $S(\alpha) = \alpha^{-1}$ (\emph{coinverse}), $S_0[D]$ becomes a Hopf algebra (cf.\  \cite{Water}, §1.4).
In particular, $G := \Spec(S_0[D])$ is a diagonalizable group scheme by \cite{Water}, §2.2, and also the Cartier dual of the constant group scheme associated to $D$ by \cite{Water}, §2.2 and §2.3. Thus we have an induced action of $G$ on  $\Spec(S)$, corresponding to the coaction of $S_0[D]$ on $S$ (also see \cite{SGA3}, I 4.7.3).

As $G$ is a diagonalizable group scheme, characters of $G$, i.e.\ homomorphisms $\chi\colon G \to \BG_m = \Spec(\CO_{S_0}[T, T^{-1}])$, correspond to invertible elements in $\CO(G) = S_0[D]$, i.e.\ elements $\alpha \in S_0[D]$, such that $\Delta(\alpha) = \alpha \otimes \alpha$ by \cite{Water}, §2.1. Elements of that type are called \emph{group-like}.

\begin{proposition}\label{prop:S_f^G=S_(f)}
For all homogeneous $f \in S$, it holds that
\begin{align*}
    (S_f)^G \equ S_{(f)},
\end{align*}
where we view both $(S_f)^G$ and $S_{(f)}$ as subrings of $S_f$.
\end{proposition}

\begin{proof}
Since the $G$-coaction sends a homogeneous $f$ to $f\otimes\chi^{\deg(f)}$, it extends to $S_f$,
so replacing $S$ by $S_f$ (where $f$ is a unit) reduces the claim to the case $f=1$. The action of $G$ on $X = \Spec(S)$ is given in terms of a morphism of schemes $\alpha\colon G \times_{S_0} X \to X$ and comes together with a projection $\pr\colon G \times_{S_0} X \to X$. The corresponding coactions are given by $\alpha^\ast\colon S \to S_0[D] \otimes_S S$ and $\pr^\ast\colon S \to S_0[D] \otimes_S S$ and give rise to the exact sequence
    \begin{align*}
        S^G \to \left( S 
\ \overset{\alpha^\ast}{\underset{\mathrm{\pr^\ast}}{\rightrightarrows}}\ 
S \otimes_{S_0} S_0[D]\right).
    \end{align*}
    Thus, we can deduce that
    \begin{align*}
        S^G &\equ \{s \in S \mid \alpha^\ast(s) \equ \pr^\ast(s)\} \\
        &\equ \{s \in h(S) \mid s \otimes \deg(s) \equ  s \otimes 1\} \\
        &\equ \{s \in h(S) \mid s \otimes \chi^{\deg(s)} \equ  s \otimes \chi^0\},
    \end{align*}
    i.e. $s \in S^G$ if and only if $\deg(s) = 0$.
\end{proof}

Using the \emph{magic of potions} (\cite{May}, Proposition 2.11), the glueing of the $\Spec(S_{(f)})$ is well-behaved and yields a scheme (\cite{May}, Construction 3.1).

\begin{definition}[Multigraded spectrum of a multigraded ring]\label{def:multihom_spec}\makebox{}{}\\
The \emph{multigraded spectrum of $S$} is defined as the scheme
\begin{align*}
    \Proj^D(S) \deq \bigcup_{f \in \Rel^D(S)} \Spec(S_{(f)}).
\end{align*}
\end{definition}

\begin{remark}
    Note that this is not the way \cite{BS} defined the multigraded Proj construction. The reason for this is that their approach only works if the group scheme $G$ can be identified with its set of $S_0$-points $G(S_0)$. However, this fails if, for example, $G$ is finite cyclic of order $n$ and $S_0$ is a field of characteristic prime to $n$, that does not contain all $n$th roots of unity\footnote{This example was given to me by an anonymous referee.}.
\end{remark}


The next example can be seen as some sort of minimal example for the peculiarities arising from this construction. 

\begin{example}\label{standard_ex_origin}
Consider $S = \BC[x, y, z]$ graded by $\BZ^2$, where $x \mapsto (1, 0)$, $y \mapsto (0, 1)$, $z \mapsto (1, 1)$ and $S_+ = (xy, xz, yz)$. The action of $G = \Spec(\BC[\BZ^2]) = \BG_m^2$ on $\Spec(S)$ is given by $(\lambda , \mu)(a, b, c) \mapsto (\lambda a, \mu b, \lambda \mu c)$. It holds that $S^G = S_0 =\BC$, thus we cannot use affine GIT (for the action on $\Spec(S)$) to construct a good quotient. Instead, we look at the action of $G$ on $\Spec(S_f)$ for relevant $f$.
Consider the affine projection morphism \begin{align*}
    \pi_+ \colon \Spec(S) \setminus V(S_+) \, \to\, \Proj^D(S) .
\end{align*}
For $U = \Spec(S_{(xz)})$, $V = \Spec(S_{(yz)})$, and $X = \Proj^D(S)$, Proposition~\ref{prop:S_f^G=S_(f)} implies for 
    \begin{align*}
        \CO_X(\pi_+^{-1}(U)) &\equ (S_{xz})^G \equ S_{(xz)} \equ \BC[\frac{xy}{z}],\ \text{ and likewise} \\
        \CO_X(\pi_+^{-1}(V)) &\equ  (S_{yz})^G \equ S_{(yz)}  \equ \BC[\frac{xy}{z}].
    \end{align*}
    Therefore, distinct open affine subsets might have the same regular functions.
    Note that for $W = \Spec(S_{(xy)})$ we get
    \begin{align*}
         \CO_X(\pi_+^{-1}(W)) \equ  S_{(xy)}  \equ \BC[\frac{z}{xy}].
    \end{align*}
    In particular, the invariant rational functions of $\Quot(S)$ are generated by $\frac{xy}{z}$ and its inverse. Hence, one can already deduce that $\Proj^D(S)$ is a $\BP^1$ with doubled origin (because the rational function $xy/z$ has two zeroes). Hence, we get a $\BP^1$ if we throw away one of $U$ or $V$. 
    
    Also note that there are small modifications of the grading that do not change Proj. For example, let $S' = \BC[x, y, z]$ with $D = \BZ^2$, where $x \mapsto 2 e_1$, $y \mapsto e_2$ and $z \mapsto e_1 + e_2$. Then still $S'_+ = (xy, xz, yz)$, so $\Proj(S) = \Proj(S')$. The reason for that is quite simple: even though the equation changes from $\deg(xy) = \deg(z)$ to $\deg(xy^2) = \deg(z^2)$, the localizations cannot see this, i.e.\ $S_{(z^2)} = S_{(z)}$ and $S_{(xy^2)} = S_{(xy)}$. 
\end{example}

\begin{remark}[$\Proj^D(S)$ is an $S_0$-scheme]\label{rem:proj_base}
$X = \Proj^D(S)$ can be identified as a relative scheme over the base scheme $\Spec(S_0)$ in a natural way: Every $f \in \Rel^D(S)$ gives rise to a canonical ring homomorphism $S_0 \to S_{(f)}$ via $s \mapsto \frac{s}{1}$. This in turn induces a morphism $\Spec(S_{(f)}) \to \Spec(S_0)$ of affine spectra. Gluing all these morphisms together, we get a morphism $\Proj^D(S) \to \Spec(S_0)$.  
\end{remark}

Open subsets of $\Proj^D(S)$, defined by homogeneous elements, are covered by basic open subsets $\Spec(S_{(f)})$ for relevant $f \in S$.

\begin{lemma}[Distinguished open subsets]\label{lem:distinguish}
Let $f \in \Rel^D(S)$ and let $h \in S_d$, $d \in D$, be homogeneous such that $hf\neq 0$.
Then 
\begin{enumerate}[label=(\roman*)]
    \item $hf \in \Rel^D(S)$.

    \item the morphism $\pi_f\colon \Spec(S_f) \to \Spec(S_{(f)})$ is open.
    
    \item $h$ defines a principal open subset $H = D(hf)$ in $\Spec(S_f)$, whose image in $\Spec(S_{(f)})$ is open.
    
    \item the image of $H$ in $\Spec(S_{(f)})$ is equal to $\Spec(S_{(hf)})$.
    
    \item the sets 
    \begin{align*}
        D^+(h) \deq \bigcup_{\substack{f \in \Rel^D(S) \\ f \in (h)}} \Spec(S_{(f)}) \equ \bigcup_{\substack{f \in \Rel^D(S) \\ h \mid f}} \Spec(S_{(f)}) \subseteq \Proj^D(S)
    \end{align*}
    are open subschemes of $\Proj^D(S)$, which obey the usual formal properties, i.e.
        \begin{enumerate}[label=(\arabic*)]
            \item $D^+(hh') \equ D^+(h) \cap D^+(h')$.
            \item $\bigcup_{h \text{ homogeneous}} D^+(h) \equ \Proj^D(S)$.
        \end{enumerate}
\end{enumerate}
\end{lemma}

\begin{proof}
\begin{enumerate}[label=(\roman*)]
    \item Obviously, $hf$ is homogeneous. Every homogeneous divisor of $f$ is also a homogeneous divisor of $hf$, hence $D^f \subseteq D^{hf}$. As $D^f$ has finite index, $D^{hf}$ also has to have finite index. 
    
    \item Let $U \subseteq \Spec(S_f)$ be open. Then $\pi_f(U) = p\restrict_{D(f)} (U)$ is open by Proposition~\ref{prop:S_f^G=S_(f)}.
    
    \item We have $\Spec(S_f) = D(f)$, so $x \in \Spec(S_f)$ corresponds to a prime ideal $\Fp \lhd S$ not containing $f$. Then $h$ defines the set $\{P \lhd S_f \mid \frac{h}{1} \not\in P\}$ corresponding to 
    \begin{align*}
        \{Q \lhd S \mid f, h \not\in Q\} \equ \{Q \lhd S \mid f \not\in Q\} \cap \{Q \lhd S \mid h \not\in Q\} ,
    \end{align*}
    where the latter is just the principal open set $D(hf) \subset \Spec(S)$. Thus $h$ defines the principal open subset $\Spec(S_{hf}) \subseteq \Spec(S_f)$. Now apply (ii). 
    
    \item  This follows immediately from (iii).
    
    \item For (1), we observe that on the left-hand side, we take the union over all relevant elements that are divisible by $h$ and $h'$, which is exactly the right-hand side. Regarding (2), we see that for $h = 1$ the set $D^+(h)$ is exactly $\Proj^D(S)$.
\end{enumerate}
\end{proof}

\section{$D$-graded Proj via quotient constructions}
Next, we want to examine the morphism of schemes $\pi_f \colon \Spec(S_f) \to \Spec(S_{(f)})$, which is induced by the ring homomorphism $S_{(f)} \to S_f$. Most importantly, we show that $\Proj^D(S)$ coincides with the geometric quotient obtained from gluing the affine geometric quotients \(\Spec(S_f) /  G\cong\Spec(S_{(f)})\).

\begin{definition}\label{def:veronese_ring}
Let $H \le D$ be a subgroup. The \emph{Veronese subring} of $S$ with respect to $H$ is defined as
\begin{align*}
    S_H \deq \bigoplus_{d \in H} S_d \subseteq S.
\end{align*}
$S_H$ is obviously graded by $H$ and $f \in S_H$ is $D$-homogeneous if and only if it is $H$-homogeneous.
\end{definition}

The following fact from \cite{BS} (page 6) states that periodic rings have \emph{enough} units to form a Laurent polynomial algebra with $\rank(D)$ variables. 

\begin{lemma}\label{lem:laurent_pol_algebra}
Let $S$ be periodic. Then there is a free abelian subgroup $F \le D$ of finite index such that the Veronese subring $S_F = \bigoplus_{d \in F} S_d$ is a \emph{Laurent polynomial algebra} $S_0[T_1^{\pm 1}, \ldots, T_r^{\pm 1}]$. 
\end{lemma}

\begin{proposition}\label{prop:index_integral}
Let $H \le D$ be of finite index. Then $S|S_H$ is integral.
\end{proposition}

\begin{proof}
It suffices to consider a homogeneous element $f\in S_d$. The number of cosets $i d H$ for $i \in \BN$ has to be finite. Therefore, there is an element $N \in \BN$ such that $N d \in H$. But 
\begin{equation*}
    \deg_D(f^N) \equ \deg_D(f) + \ldots + \deg_D(f) \equ Nd
\end{equation*}
that is, $f^N \in S_H$. Thus $f$ is integral over $S_H$.
\end{proof}

Now we can examine $\pi_f$. Note that the corresponding statement in \cite{BS} (Lemma 2.1) is stated for arbitrary rings $S$, but the proof uses GIT for $S$  being a finitely generated algebra over a field. We will work with the following more general notion (cf.\   \cite[\href{https://stacks.math.columbia.edu/tag/04AE}{Definition 04AE}]{stacks-project}). In particular, our statement is new and strictly generalizes \cite{BS}, Lemma 2.1.

\begin{definition}\label{def:geom_quot_gen}
    Let $T$ be a scheme and $B$ an algebraic space over $T$. Let $j \colon R \to U \times_B U$ be a pre-relation. A morphism $\Phi\colon U \to X$ of algebraic spaces over $B$ is called a \emph{geometric} quotient, if
    \begin{enumerate}[label=(\arabic*)]
        \item $\Phi$ is an orbit space.
        \item condition (1) holds universally, i.e. $\Phi$ is universally submersive, and
        \item the functions on $X$ are the $R$-invariant functions on $U$.
    \end{enumerate}
\end{definition}

The concept of a geometric quotient is closely related to that of a pseudo $G$-torsor, but the latter is a strictly stronger notion (cf.\  \cite[\href{https://stacks.math.columbia.edu/tag/0498}{Definition 0498}]{stacks-project}). 

\begin{definition}\label{def:pseudo_G_torsor}
    Let $D$ be a finitely generated abelian group and $S$ a $D$-graded ring. Let $G = \Spec(S_0[D])$ and denote the group action on $U_f = \Spec(S_f)$ by $\alpha_f \colon G \times_{\Spec(S_0)} U_f \to U_f$.
    \begin{enumerate}[label=(\arabic*)]
        \item  We say that $U_f$ is a \emph{pseudo $G$-torsor} under $G$, if the induced morphism of schemes 
        \begin{align*}
            \Phi_f \colon G \times_{\Spec(S_0)} U_f \to U_f \times_{\Spec(S_0)} U_f,\ (g, x) \mapsto (\alpha_f(g, x), x)
        \end{align*}
   is an isomorphism of $S_0$-schemes. The map $\Phi_f$ is often called \emph{shear map} in the literature.
   \item $U_f$ is called a \emph{trivial pseudo $G$-torsor}, if $U_f$ is a pseudo $G$-torsor such that there exists an $G$-equivariant isomorphism $G \to U_f$ over $\Spec(S_0)$, where $G$ acts by left multiplication.        
    \end{enumerate}
\end{definition}

As the map $\Phi_f$ is an isomorphism if and only if
\begin{align*}
    \varphi_f \colon S_f \otimes_{S_0} S_f \to S_f \otimes_{S_0} S_0[D],\ g \otimes h \mapsto gh \otimes \chi^{\deg(g)}
\end{align*}
is an isomorphism, we can see that we need to require some conditions on $S_f$. Taking $1 \otimes \chi^d \in S_f \otimes_{S_0} S_f$, it is immediate that the surjectivity of $\phi_f$ is equivalent to the existence of units in $(S_f)_d$ for every degree $d \in D$. This is apparently a strictly stronger condition than relevance, where units only exist in every degree $d \in D^f$.
We will also see that the pre-relation $j$ from Definition~\ref{def:geom_quot_gen} is exactly the scheme morphism $\Phi_f$.

\begin{theorem} \label{lem:periodic_quotient} 
Let $S$ be an effectively $D$-graded periodic ring.
Then the map $\pi \colon \Spec(S) \to \Spec(S_0)$ induced by $S_0 \to S$ is a geometric quotient in the category of schemes.

In particular, for all relevant $f \in S$ (where $S$ is not necessarily periodic), the corresponding statement holds for the maps $\pi_f \colon \Spec(S_f) \to \Spec(S_{(f)})$.
\end{theorem}

\begin{proof}
We start with identifying the roles in the above Definition.
Let $T = B = \Spec(S_0)$ and $\Phi = \pi_f$ (hence $U = U_f = \Spec(S_f)$ and $X = \Spec(S_{(f)})$). Then for $R := G \times_{\Spec(S_0)} U_f$ we see that $\Phi_f \colon G \times_{\Spec(S_0)}U_f \to U_f \times_{\Spec(S_0)} U_f$, given by $(g,x) \mapsto (gx, x)$ is a morphism  of schemes (induced action on $\Spec(S_f)$), and hence $j = \Phi_f$ is a pre-relation. As $j$ is exactly the morphism that corresponds to the induced action on $\Spec(S_f)$, that is, $\Phi_f$, we can apply Proposition~\ref{prop:S_f^G=S_(f)} directly to deduce (3). \\

Let $D' \le D$ denote the finite index subgroup of $D$ that is given by the degrees of homogeneous units, and denote the corresponding free group from Lemma~\ref{lem:laurent_pol_algebra} by $F$. 
We have a short exact sequence
 \begin{align*}
        0 \to F \to D \to D/F \to 0 ,
\end{align*}
giving rise to a short exact sequence of diagonalizable group schemes \begin{align*}
        0 \to G_{D/F} \to G_D \to G_F\to 0
\end{align*}
 by Proposition~\ref{prop:G_functor_represent}, where we denote $G_M = \Spec(S_0[M])$ for an abelian group $M$. It is also clear that we have ring maps $S_0 \to S_F \to S$, giving rise to a factorization of $\pi$ into $\Spec(S) \to \Spec(S_F) \to \Spec(S_0)$. 
We aim to show that:
    \begin{itemize}
        \item $\pi_1 \colon \Spec(S) \to \Spec(S_F)$ is a geometric quotient w.r.t $G_{D/F}$.
        \item $\pi_2 \colon \Spec(S_F) \to \Spec(S_0)$ is a geometric quotient w.r.t. $G_F$.
    \end{itemize}
    So let's start with $\pi_1$. By Proposition~\ref{prop:S_f^G=S_(f)} we know that $S^{G_{D/F}}$ coincides with the degree zero part of $S$ viewed as $D/F$-graded ring, which is exactly $S_F$, showing (3). Since  $S|S_F$ is integral by Proposition~\ref{prop:index_integral} and surjective by \cite[\href{https://stacks.math.columbia.edu/tag/00GQ}{Tag 00GQ}]{stacks-project}, we have an integral surjective morphism $\pi_1 \colon \Spec(S) \to \Spec(S_F)$ with a $G_{D/F}$--invariant action. 
    By \cite[\href{https://stacks.math.columbia.edu/tag/01WM}{Lemma 01WM}]{stacks-project}, $\pi_1$ is also universally closed. Then by  \cite[\href{https://stacks.math.columbia.edu/tag/0AAU}{Lemma 0AAU}]{stacks-project} $\pi_1$ is submersive. As $\pi_1$ remains surjective and closed after base change, it also remains submersive after base change, and (2) holds true. 
    Finally, regarding (1), we apply \cite[\href{https://stacks.math.columbia.edu/tag/049Z}{Lemma 049Z}]{stacks-project}, having the conditions (i)-(iii). Condition (iii) holds since $\pi_1$ is surjective, and condition (i) holds since
\[
S^{G_{D/F}}=S_F.
\]
For condition (ii), the equality $S^{G_{D/F}}=S_F$ implies that the fibers over geometric points are precisely the $G_{D/F}$-orbits.
Hence, the induced morphism
\[
j_F\colon
G_{D/F}\times_{\operatorname{Spec}(S_0)}\operatorname{Spec}(S)
\longrightarrow
\operatorname{Spec}(S)\times_{\operatorname{Spec}(S_F)}\operatorname{Spec}(S)
\]
is surjective. Thus (1) holds. \\

Regarding $\pi_2$, it holds that $S_F \cong S_0[F]$, as $S_F$ is a $F$-graded Laurent polynomial ring and is of finite type over $S_0$ (since $\rank(D)$ is finite). Furthermore, the action of $G_F$ on $\Spec(S_F)$ is free and transitive, i.e.\ $\Spec(S_F) = G_F$ (and the corresponding shear map from Definition~\ref{def:pseudo_G_torsor} is an isomorphism). Thus $\pi_2$ is a Zariski locally trivial (pseudo) $G_F$-torsor and hence the projection $\Spec(S_F) \to \Spec(S_0)$ is a geometric quotient. 
As the composition of geometric quotients is a geometric quotient, $\pi\colon \Spec(S) \to \Spec(S_F) \to \Spec(S_0)$ is a geometric quotient by the action of $G_D$ following the exact sequence above. Since both $\pi_1$ and $\pi_2$ are affine, the geometric quotient actually is a scheme by construction (cf.\   \cite[\href{https://stacks.math.columbia.edu/tag/03JH}{Lemma 03JH}]{stacks-project}).
\end{proof}

\begin{remark}
An anonymous referee pointed out a related approach based on \cite{DG}, Chapter III, §2, 6.1, Corollaire. 
For a finite-index
free subgroup $F\subset D$, write $E=D/F$ and let
\[
\varphi\colon \operatorname{Spec}(S)\to G_F
\]
be the morphism determined by homogeneous units whose degrees generate
$F$. Setting
\[
Y\deq \varphi^{-1}(e_F),
\]
where $e_F\colon \operatorname{Spec}(S_0)\to G_F$ denotes the neutral
section, the referee suggests applying Corollaire~6.1 to the induced
action of the finite locally free diagonalizable group scheme
$G_E$ on $Y$. This gives a quotient for the
$G_{D/F}$-action on $Y$. Together
with a $G_D$-equivariant identification of the form
\[
\operatorname{Spec}(S)\simeq
(G_D\times Y)/G_E,
\]
this may lead to a description of the quotient for the $G_D$-action.
We do not pursue this approach here, since the argument above gives
the desired geometric quotient directly in the sense of
Definition~\ref{def:geom_quot_gen}.
\end{remark}

For the pseudo $G$-torsor structure, we need stronger assumptions on $S_f$.

\begin{theorem}\label{thm:quotient_pi_f_general}
    Let $S$ be a $D$-graded ring, $f \in \Rel^D(S)$ and $\pi_f\colon \Spec(S_f) \to \Spec(S_{(f)})$ be the map induced by $S_{(f)} \to S_f$. Then $\Spec(S_f)$ is a pseudo $G$-torsor if and only if $\deg( (S_f)^\times) = D$ (or equivalently $D^f = D$) and $S$ is integral.
    In particular, if $\deg( (S_f)^\times) = D$ holds and $S$ is integral, then $\pi_f$ is Zariski locally trivial for the base change $S_0 \to S_{(f)}$.
\end{theorem}

\begin{proof}
First note that the new condition $\deg( (S_f)^\times) = D$ ensures that we have a unit $u_d \in (S_f)^\times$ of degree $d$ for every $d \in D$. Hence $\supp( (S_f)^\times) =D$ and $1 \in S_d S_{-d}$ for all $d \in D$, i.e. $S_f$ is \emph{strongly graded} (cf.\  \cite{Dade80}, Proposition 1.6).
Let $\deg( (S_f)^\times) = D$ hold, where $S$ is integral.
Regarding injectivity, since $S_f$ is graded, the kernel of $\varphi_f$ is graded as well. Thus, it is enough to check homogeneous elementary tensors.
An element $g \otimes h$ lies in the kernel of $\varphi_f$ if and only if $gh = 0$, since $S_f$ is graded. As we assume $S$ to be integral, we can deduce that $g=0$ or $h=0$. Now as $0 \otimes 0 = a \otimes 0 = 0 \otimes b$, we see that $g \otimes h \in \ker(\varphi_f)$ if and only if $g \otimes h$ is the $0$-tensor. Hence $\varphi_f$ is injective. 
Conversely, the injectivity of $\varphi_f$ forces $S$ to be an integral domain.

For surjectivity, we have to show that an arbitrary element $g \otimes \chi^d \in S_f \otimes_{S_0} S_0[D]$ lies in the image of $\varphi_f$. By assumption, there exists a unit $u_d \in (S_f)_d$, hence for $g = u$ and $h = u^{-1} g$ we see that $\varphi_f(gh) = g \otimes \chi^d$. 
If $\varphi_f$ is surjective, this forces units in $S_f$ to exist in every degree $d\in D$.

For the local triviality, we have to show that there is an isomorphism $\varphi\colon S_f \to S_{(f)}[D]$. First, by assumption, we can choose homogeneous units $u_d \in (S_f)_d$ for every degree $d \in D$. We define $\varphi$ by $q \mapsto \frac{q}{u_{\deg(q)}} \chi^{\deg(q)}$. For $q \chi^d \in S_{(f)}[D]$, we see that $\varphi(u_d q) = q \chi^d$ and $\varphi$ is surjective.
Regarding injectivity, let $q = \sum_{d\in D} q_d$ such that $\varphi(q) = 0$, i.e.\ $\sum_{d\in D} \frac{q_d}{u_d} \chi^d=0$. 
As $\chi^d \neq 0$ and the various $\chi^d$ are linearly independent by Proposition~\ref{prop:group_like} (since $S$ is integral), we deduce that $\frac{q_d}{u_d} = 0$ for all $d \in D$ and hence $q_d = 0$ for all $d$ and $\varphi$ is injective. Finally, for $g \in G$ and $d := \deg(q)$ it holds that
    \begin{align*}
        \varphi(g \bullet q) \equ  \varphi(\chi^d(g) q) \equ
        \chi^{d}(g) \varphi(q) \equ
        \chi^{d}(g) \frac{q}{u_{d}} \chi^{d} \equ g \odot (\frac{q}{u_{d}} \chi^{d}) \equ g \odot \varphi(q),
    \end{align*}
    where $\bullet$ denotes the induced $G$-action on $S_f$ and $\odot$ denotes the induced $G$-action on $S_{(f)}[D]$ (that is induced via the base change $S_0 \to S_{(f)}$). Hence $\varphi$ is also $G$-equivariant.
\end{proof}

The above discussion gives rise to the following definition.

\begin{definition}
    Let $D$ be a finitely generated abelian group and $S$ a $D$-graded ring. We call a relevant element $f \in \Rel^D(S)$ \emph{strongly relevant}, if $S_f$ contains a homogeneous unit in every degree $d \in D$, i.e., when $D^f = D$.
\end{definition}

We can interpret the shear map in terms of properties of the action of $G$ on $U_f$.

\begin{lemma}
$\Phi_f$ is injective if and only if the action of $G$ on $\Spec(S_f)$ is free.
Moreover, $\Phi_f$ is surjective if and only if the action of $G$ on $\Spec(S_f)$ is transitive.
\end{lemma}

\begin{proof}
    If the action is free and there exist $x, y \in U_f$ and $g, h \in G$ such that $(gx,x) = (hy,y)$, clearly $x= y$ and hence $gx = hx$. Therefore, $g = h$ by freeness. Conversely, let $\Phi_f$ be injective and let $g \in G, x \in U_f$ such that $gx = x$, i.e.\ $(gx, x) = (1_G x, x)$. By injectivity, we deduce that $g = e$ and thus, the action is free.

    For the second statement, it holds that $(x, y) \in \im(\Phi_f)$ if and only if there exists a $g \in G$ such that $gy = x$, i.e.\ such that $(x, y) = (gy, y)$.
\end{proof}

For relevant $f \in S$, we define $G^f := \Spec(S_0[D^f])$.

\begin{corollary}\label{cor:G^f-torsor}
    Let $D$ be a finitely generated abelian group, $S$ a $D$-graded integral domain 
    and $f \in S$ relevant. Then $\pi_f \colon \Spec(S_f) \to \Spec(S_{(f)})$ is a $G^f$-torsor.
    In particular, if $\supp(S_f) = D^f$, then $\pi_f$ is Zariski locally trivial for the base change $S_0 \to S_{(f)}$.

Hence, for an integral $D$-graded domain $S$ and relevant $f \in S$, $\pi_f$ is a pseudo $G$-torsor if and only if $\supp((S_f)) = D^f = D$.
\end{corollary}

\begin{proof}
    We have to show that the map
    \begin{align*}
    \varphi_f\colon    S_f \otimes_{S_0} S_f \to S_f \otimes_{S_0} S_0[D^f],\ g \otimes h \mapsto gh \otimes \chi^{\deg(g)}
    \end{align*}
    is an isomorphism. However, as we assume the grading to be effective, we know that $\varphi_f$ is an isomorphism if and only if $\deg((S_f)^\times) = D^f$ by Theorem~\ref{thm:quotient_pi_f_general}. But this holds by definition of $D^f$, as every divisor of some power of $f$ corresponds to a unit in $S_f$. Alternatively, this statement is equivalent to the second part of the proof of Theorem~\ref{lem:periodic_quotient} and hence true.

    Now consider the base change $S_0 \to S_{(f)}$. Then $\pi_f$ is a trivial $G^f$-torsor if and only if $S_f \otimes_{S_{(f)}} S_{(f)} \cong S_{(f)} \otimes_{S_0} S_0[D^f] = S_{(f)}[D^f]$ as $S_{(f)}$-algebras. For all $d \in D^f$, we choose a unit $u_d \in ((S_f)_d)^\times$ and consider the following morphism of $S_{(f)}$-algebras
    \begin{align*}
        \varphi\colon S_f \to S_{(f)}[D^f],\ q \chi^d \mapsto \frac{q}{u_{\deg(q)}} \chi^d,
    \end{align*}
    which is well-defined by assumption (i.e.\ $\supp(S_f) = D^f$).
    Let $q \cdot \chi^d \in S_{(f)}[D^f]$. Then
    \begin{align*}
        \varphi(u_d q) \equ \frac{q u_d}{u_d} \chi^d \equ q \chi^d,
    \end{align*}
    i.e.\ it is surjective. Regarding injectivity, let $q = \sum_{d\in D^f}$ such that $\varphi(q) = 0$. Then $\sum_{d\in D^f} \frac{q_d}{u_d} \chi^d=0$. As $\chi^d \neq 0$ and the various $\chi^d$ are linearly independent (by Proposition~\ref{prop:group_like}), we deduce that $\frac{q_d}{u_d} = 0$ for all $d \in D^f$ and hence $q_d = 0$ for all $d$ and $\varphi$ is injective. Finally, $\varphi$ is $G^f$-equivariant by the same reasoning as in Theorem~\ref{thm:quotient_pi_f_general}.
\end{proof}

\begin{example}\label{ex:double_origin_3rd}
    The above conditions seem to be satisfied in many cases. For example, consider again the ring $S = \BC[x, y, z]$ graded by $\BZ^2$, where $x \mapsto (1, 0)$, $y \mapsto (0, 1)$ and $z \mapsto (1, 1)$. It holds $S_+ = (xy, xz, yz)$, and we have already shown that every generator is strongly relevant in Example~\ref{ex:double_origin_P^1_2nd}, so that $\pi_f$ is locally trivial in the Zariski topology for $f = xy, xz, yz$.
    In particular, as $D = \BZ^2$ has no torsion, we conclude that $\Proj^D(S)$ is a pseudo $G$-torsor. 
    Also note that Example~\ref{ex:standard_ex} (1) satisfies those conditions, whereas (2) does not.
\end{example}

\begin{remark}
    Note that our notion of strong relevance coincides with the notion of maximal relevance in \cite{May}, Section 5.4. Also note that their Proposition 5.13 is close to Corollary~\ref{cor:G^f-torsor}, with one big difference. We work over the base $\Spec(S_0)$, while \cite{May} show that $\pi_f$ is a trivial pseudo $G$-torsor over the base $\Spec(S_{(f)})$. Over the latter base, no integrality condition on $S$ is needed.
\end{remark}

The following definition is due to \cite{BS}.

\begin{definition}\label{def:irrelevant_ideal}
We call $V(S_+) \subseteq \Spec(S)$ the \emph{irrelevant subscheme} and $D(S_+) \subseteq \Spec(S)$ the \emph{relevant locus}. We denote the affine projection $\Spec(S) \setminus V(S_+) \to \Proj^D(S)$ by $\pi_+$.
\end{definition}

\begin{proposition}\label{prop:irrel_quotient}
Let $S$ be a multigraded ring, graded by $D$ with free abelian part $F$.
\begin{enumerate}[label=(\alph*)]
    \item Let $S$ be a finitely generated $S_F$-module.
    Then the affine projection morphism
\begin{align*}
    \pi_+ \colon \Spec(S) \setminus V(S_+) \, \to\, \Proj^D(S)
\end{align*}
is a geometric quotient for the induced action. 

    \item Let $S$ be an integral domain 
    and $\deg((S_f)^\times) = D$ holds for all relevant $f \in S$. Then the affine projection morphism
\begin{align*}
    \pi_+ \colon \Spec(S) \setminus V(S_+) \, \to\, \Proj^D(S)
\end{align*}
is a pseudo $G$-torsor that is Zariski locally trivial over the base $\Spec(S_{(f)})$.
\end{enumerate}
\end{proposition}

\begin{proof}
    The statements follow from Theorem~\ref{lem:periodic_quotient} and Theorem~\ref{thm:quotient_pi_f_general} and gluing, as the conditions are local in either case.
\end{proof}

The affine projection comes with the expected properties. The following statement is the generalization of a standard result from affine GIT (cf.\  \cite{Hos}, Corollary 4.31).

\begin{proposition}\label{prop:orbit_closed}
Consider $\pi_+\colon D(S_+) \to \Proj^D(S)$. Then
    \begin{align*}
        \pi_+(x) = \pi_+(y) \ \iff \ \overline{G \cdot x} \cap \overline{G \cdot y} \neq \emptyset 
    \end{align*}
\end{proposition}

\begin{proof}
\Ra As $\pi_+$ is a geometric quotient by Theorem~\ref{lem:periodic_quotient}, fibers of $\pi_+$ are exactly orbits of $G$, so $\pi_+(x) = \pi_+(y)$ yields that $x$ and $y$ lie in the same orbit. Hence $\overline{G \cdot x} = \overline{G \cdot y}$ and thus their intersection cannot be empty (since $x, y$ are contained). 

\La 
Let $z\in\overline{G\cdot x}\cap\overline{G\cdot y}$. Since $z\in D(S_+)$, there exists a relevant $f\in S$ such that
$z\in\Spec(S_f)$. As $\Spec(S_f)$ is open, we can choose
\[
x'\in G\cdot x\cap\Spec(S_f),
 \text{ and }
y'\in G\cdot y\cap\Spec(S_f).
\]
Assume that $\pi_f(x')\neq\pi_f(y')$. Then there exists a
$G$-invariant regular function $h\in S_{(f)}$ such that
\[
h(\pi_f(x'))\neq h(\pi_f(y')).
\]
Since $h$ is $G$-invariant, it is constant on $G$-orbits and hence on
their closures. Thus
\[
h(\pi_f(x'))=h(\pi_f(z))=h(\pi_f(y')),
\]
a contradiction.
\end{proof}


\begin{example}[Projective Space]\label{exa:special_proj}
Let $S = \BC[x_0, \ldots, x_n]$ be a $\BZ$-graded polynomial ring such that $\deg_D(x_i) =1$. Hence, $S_+ = (x_0, \ldots, x_n)$, therefore $\Proj^D(S)$ is the classical Proj, which we denote by $\Proj^\BN(S)$. The group action by $G = \Spec(S_0[D])= \BG_m$ on $\Spec(S) = \BA^{n+1}$ is given by 
\begin{align*}
    (\lambda, (a_0, \ldots, a_n)) \mapsto (\lambda a_0, \ldots, \lambda a_n).
\end{align*}
In particular, this coincides with the classical Proj construction. Also note that each weighted projective space is naturally given by a $D$-graded Proj. If $X = \BP(q_0, \ldots, q_n)$ is a weighted projective space with weights $q_i$, then $S = \BC[x_0, \ldots, x_n]$ with $\deg(x_i) = q_i$ (cf.\  \cite{Cox}, page 3). In particular $X = \Proj^\BN(S)$.
\end{example}

There are surprisingly many gradings yielding a trivial $D$-graded Projectivization (i.e.\ $|\Proj^D(S)| = 1$).

\begin{proposition}\label{prop:Proj_trivial_r=n}       
Let $S = \BC[x_1, \ldots, x_n]$ be $D$-graded such that the degrees of the $x_j$ are pairwise linearly independent and $r:=\rank(D) = n$. Then $\Proj^D(S) = \{(0)\}$.
\end{proposition}

\begin{proof}
Since all the $\deg(x_j)$ are pairwise linearly independent and $\rank(D) = n$, we conclude that relevant elements must have exactly $n$ linearly independent divisors. Hence $S_+ = (x_1 \cdot \ldots \cdot x_n)$, i.e.\ there is only one affine open giving Proj. It holds $S_{(x_1 \cdot \ldots \cdot x_n)} = \BC$ and hence $\Proj^D(S) = \Spec(\BC) = \{(0)\}$. 
\end{proof}

The following examples give some well-known spaces.

\begin{example}\label{ex:P1_first}
    \begin{enumerate}[label=(\arabic*)]
    \item For $S = \BC[x_0, \ldots, x_{n+1}]$ and $D = \BZ^2$ such that $\deg(x_0) = e_2$ and $\deg(x_i) = e_1$ for $i=1, \ldots, n+1$, we deduce that $\Proj^D(S) = \BP^n$. Again, all relevant elements are strongly relevant.

    \item Products of projective spaces can be easily seen as $D$-graded Proj schemes. If $X = \BP^n\times\BP^m$, then for $S = \BC[x_0, \ldots, x_n; y_0, \ldots, y_m]$ and $\deg(x_i) = e_1$ and $\deg(y_j) = e_2$ we see that $\Proj^{\BN\times\BN}(S) = \BP^n\times\BP^m$ (cf.\  \cite{Cox}, page 3).
    \end{enumerate}
\end{example}

We can also show $\Proj^D(S) = (0)$ in Proposition~\ref{prop:Proj_trivial_r=n} by computing $\dim(\Proj^D(S))$. Note that this fact was already stated in \cite{KU}, Lemma 3.6 (2), but the proof given there is not sufficient.

\begin{lemma}\label{lem:dim_proj}
Let $S$ be a noetherian effectively $D$-graded ring (cf.\ Definition~\ref{def:grading_effective}), such that for all relevant $f \in S$, we have $(\bigcup_{m\in \BN}\Ann_S(f^m))\cap S_F = \{0\}$. Then it holds
\begin{align*}
\dim(\Proj^D(S)) \equ \dim(\Spec(S)) - \rank(D)
\end{align*}
In particular, this statement holds if $S$ is integral, or when $\Ann(S_f) = \{0\}$.
\end{lemma}

\begin{proof}
Let $D =F \oplus T$ be the decomposition into free and torsion part of $D$, where $F \cong \BZ^{\rank(D)}$. 
Applying Lemma~\ref{lem:laurent_pol_algebra} to the periodic ring $S_f$ for relevant $f \in S$ immediately gives that $(S_f)_F\ \cong S_{(f)}[F]$. As $S$ is noetherian, $\dim(S_f)$ and $\dim(S_{(f)})$ are finite and we can deduce that $\dim((S_f)_F) = \dim(S_{(f)}) + \rank(D)$. 
By Proposition~\ref{prop:index_integral}, $S_f$ is integral over $(S_f)_F$ and hence $\dim(S_f) = \dim( (S_f)_F)) = \dim(S_{(f)}) + \rank(D)$ by \cite{Bosch}, Proposition 3.3/6, as the inclusion $(S_f)_F \to S_f$ is injective (by the annihilator condition). 
Since $\Proj^D(S)$ is given by the union of all $\Spec(S_{(f)}) = \Spec(S_f) \sslash G$, the claim follows.
\end{proof}


\section*{Appendix A. Details on the group scheme \texorpdfstring{$G=\Spec(S_0[D])$}{G = Spec(S0[D])}}
Most of the material in this section is well known to experts. 
However, since there is no reference treating this exact setting, we include it here for completeness.
The following statement generalizes \cite{Milne}, Lemma 12.4.

\begin{proposition}\label{prop:group_like}
Let $\Spec(S_0)$ be connected (so that $S_0$ only has the trivial idempotents $0_S$ and $1_S$).
Then the group-like elements of $S_0[D]$ are linearly independent and exactly the elements of $D$.
\end{proposition}

\begin{proof}
    As elements of $D$ are group-like, we only have to take care of the converse direction. Thus let $\alpha = \sum_i a_i d_i \in S_0[D]$ be group like, where $a_i \in S_0$ and $d_i \in D$. For $\alpha \neq 0$ it holds
    \begin{align*}
        \sum_i a_i d_i \otimes d_i \equ \sum_i a_i \Delta(d_i) \equ \Delta(\alpha) \equ \alpha \otimes \alpha \equ \sum_{i, j} a_i a_j d_i \otimes d_j
    \end{align*}
    and the $d_i \otimes d_j$ are linearly independent for all $i, j$, this implies that $a_i a_j = 0$ for $i \neq j$ and $a_i^2 = a_i$. As $S_0$ does not contain any non-trivial idempotent element, it follows that $a_j = 0$ for all $j \neq i$, that is, $\alpha = d_i \in D$ (so that $a_i = 1$ for exactly one $i$).
\end{proof}

Hence, we can interpret $D$ as the character group of the group algebra $S_0[D]$ (if $\Spec(S_0)$ is connected). In general, elements $d \in D$ correspond to characters $\chi^d: G \to \BG_m = \Spec(\CO_{S_0}[t, t^{-1}])$, so that we may view the coaction as
\begin{align*}
    S \to S \otimes_{S_0} S_0[D], \ S_d \ni f \mapsto f \otimes \chi^d, 
\end{align*}
and therefore
\begin{align*}
    S_0[D] \equ \bigoplus_{d\in D} S_0 \chi^d .
\end{align*}
In particular, we can deduce the following equation for homogeneous elements (cf.\  \cite{CRB}, Construction 1.2.2.2):
\begin{align*}
    f \in S_d \ \iff \ f(gx) \equ \chi^d(g) f(x) \ \ \text{ for all } g \in G, x \in  \Spec(S).
\end{align*}

The following result is due to \cite{Water} (Theorems on pages 5, 9, and 14 and § 2.2). 

\begin{proposition}\label{prop:G_functor_represent}
    Let $S$ be a $D$-graded integral domain, where $D$ is a finitely generated abelian group.
    \begin{enumerate}[label=(\roman*)]
        \item The affine group scheme $G = \Spec(S_0[D])$ represents the exact functor
        \begin{align*}
            G_D\colon \mathrm{Alg}_{S_0} \to \mathrm{Grp},\ \ R \mapsto \Hom_{\mathrm{Grp}}(D, R^\times).
        \end{align*}
        \item $G$ is a finite product of copies of 
        \begin{align*}
            \BG_{m,S_0} \deq \Spec(\CO_{S_0}[T, T^{-1}])
        \end{align*}
        and various 
        \begin{align*}
            \mu_{n_i, S_0} \deq \Spec(\CO_{S_0}[T]/(T^{n_i}-1)).
        \end{align*}
    \end{enumerate}
\end{proposition}

\end{document}